\documentclass{article}

\usepackage{amsmath,amssymb}
\usepackage[margin=3cm]{geometry}
\usepackage{mathtools}
\usepackage{cleveref}
\usepackage[final]{showlabels}
\usepackage[table]{xcolor}
\usepackage{graphicx}
\usepackage{tikz}
\usetikzlibrary{intersections,calc,arrows.meta}
\usepackage{pgfplots}
\pgfplotsset{compat=1.16}
\usepackage{authblk}
\usepackage{verbatim}
\usepackage{booktabs} 
\usepackage{caption}
\usepackage[subrefformat=parens]{subcaption}

\newcommand{\cp}{x} 
\newcommand{\func}{g} 
\newcommand{\con}{m}
\newcommand{\sgn}{\operatorname{sgn}}
\newcommand{\DeltaT}{\Delta t}
\newcommand{\candidate}[1]{x_{#1}}
\newcommand\movements{s}

\newcommand{\dia}[1]{D_{#1}}
\newcommand\movement[2]{h_{#1,#2}}

\newcommand{\commentinline}[2][_]{\textcolor{red}%
{%
\ifx#1_[
\else[#1: \fi
#2]}
}

\title{Explicitly Multimodal Benchmarks for Multi-Objective Optimization}
\author[1]{Ryosuke Ota}
\author[1]{Reiya Hagiwara}
\author[2]{Naoki Hamada}
\author[1]{Likun Liu}
\author[3]{Takahiro Yamamoto}
\author[1]{Daisuke Sakurai}
\affil[1]{Kyushu University}
\affil[2]{KLab, Inc.}
\affil[3]{Tokyo Gakugei University}

\date{\today}

\begin{document}


\maketitle


\begin{abstract}
    In multi-objective optimization, designing \emph{good} benchmark problems is an important issue for improving solvers.
    Controlling the global location of Pareto optima in existing benchmark problems has been problematic, and it is even more difficult when the design space is high-dimensional since visualization is extremely challenging. 
    As a benchmarking with explicit local Pareto fronts, we introduce a benchmarking based on basin connectivity (3BC) by using basins of attraction. 
    The 3BC allows for the specification of a multimodal landscape through a kind of topological analysis called the basin graph, effectively generating optimization problems from this graph. 
    Various known indicators measure the performance of a solver in searching global Pareto optima, but using 3BC can make us localize them for each local Pareto front by restricting it to its basin.  
    3BC's mathematical formulation ensures the accurate representation of the specified optimization landscape, guaranteeing the existence of intended local and global Pareto optima.
\end{abstract}


\section{Introduction}

Benchmark problems for multi-objective optimization help evaluate and understand its solvers.
A good benchmark problem would reveal solver behaviors in different situations, hinting their performance. 
In the literature of evolutionary computation, for example, we can find a number of famous benchmark suites. 
The property of multi-objective optimization problems we mainly consider is multimodality. 
Although a clear consensus on what constitutes multimodality is hardly reached, a key characteristic is the basins of attraction that trap solvers, making them fail to explore global optima \cite{mersmann2011exploratory}. (See \cref{sec: multimodality model} for more.)

In conventional multimodal multi-modal benchmarks, it is difficult to arrange the number and location of local Pareto fronts. 
Commonly used benchmark problems in \cref{sc: other existing benchmark problems} use similar constructions as used in Deb's toolkit~\cite{35_deb1999multi}.
This method, which uses several functions combined to form the objective function, may control the Pareto optima locally but cannot globally. 
That is, it may be possible to make a particular point Pareto optimum by modifying the constituting functions, but then an extra Pareto optimum will arise or disappear at another location. 
Hence, in order to find reference points, it is a standard way to run algorithms to numerous generations with a large population size, which is imprecise. 
Although there have been benchmark problems with this goal, specifying the local Pareto fronts, while avoiding unintended ones, remains a significant challenge.
When the design space becomes high dimensional, the problem intensifies since visualization becomes extraordinarily challenging. 

In this article, we propose a novel benchmark problem named the Benchmarking Based on Basin Connectivity (3BC), in which the user can completely specify the number and location of local Pareto fronts, by using basins of attractions \cite{preuss2015multimodal}. 
Therefore, the reference points for all local Pareto fronts can be analytically obtained and it makes it possible to compute performance indicators for each local Pareto front following a rigorous formulation. 
See \cref{sec: basin-wise IGDX} (\cref{fig:basin-wiseIgdx} in particular) for our demonstration.

The 3BC allows the benchmark designer to specify the multimodal landscape, which we model employing the topological analysis of preimage $f^{-1}$.
The designer can specify how the basins of attractions in the landscape connect to each other using a graph we name the \emph{basin graph} (the \emph{reachability graph} \cite{Bormann2020} is a similar yet more complex variant of our basin graph).
The 3BC generates an instance of multi-objective optimization problem from the basin graph.
As a first step in our research direction, the 3BC considers 2 objective functions. However, the authors believe the approach is generalizable to benchmarking with more objective functions.

While various indicators can measure the performance of solvers, they mainly focus only on the global Pareto optima since the local Pareto optima are generally unknown in practice.
In the paper, however, we can define some basin-wise versions since the local Pareto optima are specified in the 3BC and we can use the basin of each local Pareto optimum. 
In other words, by computing indicators only in the basin for a local Pareto optimum, we can measure the ability of solvers to search for that local Pareto optimum. 

Before this work, multimodality has been inferred usually by clustering the population of evolutionary algorithms. 
One major problem here, however, is that the induced multimodality has ambiguity: it can vary significantly depending on parameters, including the initial condition, and even random numbers picked for the clustering computation. 

Due to its mathematical formulation, the 3BC guarantees the modeled landscape to match the output.
That is, the specified local and global Pareto optima exist, forming specified landscapes, and there is no other optima.


\section{Related work}

\subsection{Multimodality model}
\label{sec: multimodality model}

As mentioned, we regard the key property of multimodality to be the difficulty for a greedy algorithm to find global optima.
Finding a satisfactory formal definition, however, has been challenging.

One common definition of multimodality is that a multi-objective optimization problem has local Pareto optima that are not global optima. Li et al.~\cite{li2023multimodal} provide a nice review.
The problem with this definition is that a multimodal problem in this sense can be solved with a greedy algorithm 
at times, which is contrary to the ideal of multimodality. 

\begin{figure}[tbh]
    \centering
    \includegraphics[width=0.3\columnwidth]{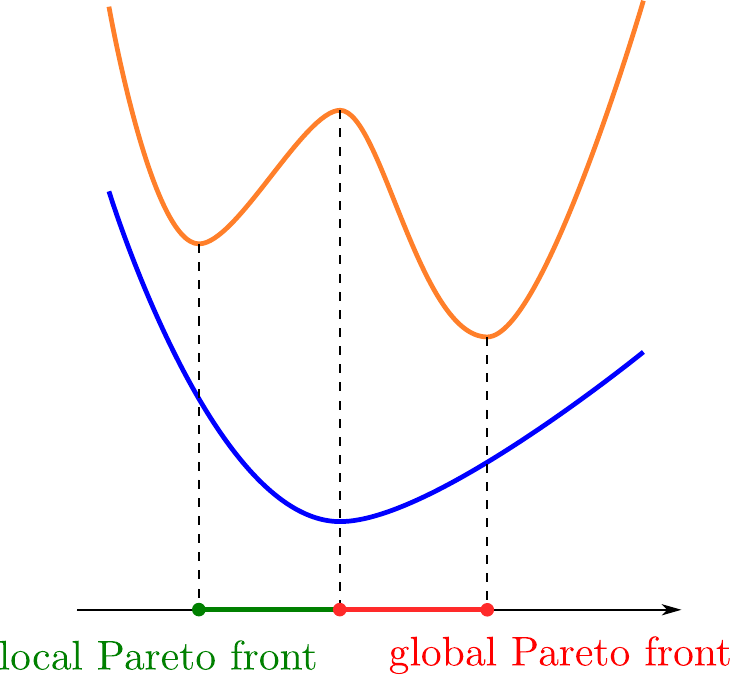}
    \caption{A 2-objetive function which has local optima but can be solved by a greedy algorithm}
    \label{fig: difference between two definitions of multimodality folklore and Huband}
\end{figure}

Huband et al.~\cite{huband2006review} considered the multimodality of each objective function.
However, there exist optimization problems that satisfies this condition and can still be solved in a greedy manner.
\Cref{fig: difference between two definitions of multimodality folklore and Huband} gives such an example.
Here, the orange objective function is by itself a multimodal objective with two pits.
However, the Pareto optimization, that also optimizes the blue objective function, gives a Pareto set that consists of the green region, local optima, and the red region, the global ones.
In such a situation, a greedy solver that follows the local optima can easily reach the global optima, which is contrary to the concept of multimodality.

According to Kerschke et al.~\cite{kerschke2016towards}, multimodality is to have multiple connected components of the local Pareto sets. 
In the example shown in \cref{fig: difference between two definitions of multimodality folklore and Huband}, the orange and blue lines are the graphs of two functions consisting of the two-objective function of one valuable. 
Each point in the green area is a local (but not global) Pareto point but the local Pareto set consists of one connected component. 

This is the base of our multimodality model.
Notice, however, that the multimodality as defined by Kerschke et al.~\cite{kerschke2016towards} can be discussed only when we are talking about points in Pareto optima, not in the whole design space.
We thus extend it so that any point in the design space can be given a modality.
We do so in the form of basins of attraction.

\subsection{Multimodal benchmark problems}

Li et al.~\cite{li2023multimodal} lists several multimodal multi-objective problems. 
Many of them have a small number of decision valuables, such as 
Omni-test~\cite{35_deb1999multi}, 
SYM-PART~\cite{36_rudolph2007capabilities}, 
TWO-ON-ONE~\cite{37_preuss2006pareto}, and 
Polygon problems~\cite{23_ishibuchi2011many}. 
In contrast, The 3BC is scalable, that is, our method can design multimodal problems with any dimension of the design space. 
Some problems proposed in Yue et al.~\cite{26_yue2019novel}, Multi-Polygon~\cite{24_ishibuchi2019scalable}, MMMOP~\cite{27_liu2018multimodal}, and SMMOP~\cite{28_tian2020multipopulation} are scalable but our method is still valid because the 3BC can also have multiple local optima yet it has more precise multimodality. 
That is, the local Pareto sets of the 3BC are completely calculated, and furthermore, we can specify the connectivity between each basin as the landscape. 
Kerschke et al.~\cite{kerschke2016towards} presented mixed sphere problems having multimodality as described in \cref{sec: multimodality model}. 

\subsection{Other existing benchmark problems}
\label{sc: other existing benchmark problems}
There are well-known benchmark problems used to evaluate the performance of evolutionary multi-objective optimization (EMO) algorithms. The ZDT test suite~\cite{zitzler2000comparison} and DTLZ test suite
\cite{deb2005scalable}
are widely used test problems. The Pareto set can be analytically obtained. 
The present benchmark of ours has the advantage that one can specify the local (and global) Pareto set more flexibly.
In the ZDT and DTLZ, one could only change the Pareto set along the direction of the distance variables.
The position variables parameterize the Pareto set.
The numbers of distance and position variables could not be changed. 

Huband et al.~\cite{huband2006review} proposed the WFG test suite to improve the ZDT and DTLZ test problems. 
The number of distance and position variables can be changed.
The WFG possesses non-separability \cite{huband2006review}, and its Pareto front is known. 
The limitation of the WFG, however, is the difficulty in tuning multiple properties of the benchmark and in controlling the Pareto set. These are interdependent, and the Pareto set cannot be explicitly defined. 
For example, the test problem WFG3, which was expected to have the degenerate Pareto front, actually was found to have a non-degenerate part in the Pareto front \cite{ishibuchi2015pareto}.
An inverted version of DTLZ and WFG has been proposed \cite{ishibuchi2017,jain2013}, yet the problem remains.
The 3BC can be seen as an improvement to the WFG that offers an explicit control of characteristics and Pareto sets simultaneously and independently.
With the 3BC we can create non-separable problems, and problems with various shapes of the Pareto front.

There are other commonly used test problems such as MaF~\cite{cheng2017benchmark}, LSMOP~\cite{cheng2016test}, and UF~\cite{zhang2008multiobjective}. 
MaF offers test problems with shapes of the Pareto front which do not appear in ZDT and DTLZ, and WFG. 
The LSMOP test problems are designed to evaluate the performance of large-scale EMO algorithms and the UF test suite has the advantage that their Pareto set and Pareto front are computed. 

Other trends of standard benchmarking include COCO bbob-biobj~\cite{hansen2021coco}, which contains 55 bi-objective problems, and a collection of real-world problems for benchmarking proposed by Tanabe et al.~\cite{tanabe2020easy}.

Our method can also make large-scale test problems with various shapes of Pareto fronts and, unlike all the above problems, has user-defined local Pareto sets and fronts. 

\section{Our multimodality: basin configuration}

We model the multimodality suitable to our motivation. 
For this, we consider the \emph{basins of attractions} \cite{preuss2015multimodal}.

Any point in the design space may reach multiple local Pareto sets without any component of the objective function increasing. 
Suppose, for example, $x$ is a point in design space and it can reach two different components $S_1$ and $S_2$ of the local Pareto set by continually improving the values. 
Then some point in $S_1$ can arrive at another point in $S_2$ by increasing the value only at the point $x$. 
From this point of view, we can consider two components $S_1$ and $S_2$ connect each other. 
In the present paper, we would like to construct a benchmark problem corresponding to the given connectivity of components of the local Pareto sets. 
We consider this type of connection between the components by using some graph called the \textit{basin graph}. 
For more details, see \cref{sec: basin graph} and \cref{sec:inputtree}. 

\section{Characteristics of our benchmark}
Our benchmark problem has characteristics that are guaranteed by the design and mathematical formulation.
They can be configured by supplying a graph structure called the basin graph as input.

\subsection{Overview} \label{sec:characteristicsoverview}

\paragraph{Graph input} The benchmark problem is generated from the user-supplied basin graph that encodes information such as the positions of local Pareto efficient points and the connectivity of the basins of attraction.
First, we consider a diamond. 
\paragraph{Explicit positioning of local and global Pareto optima}
We allow the user to specify the positions of local and global Pareto optima.
We provide the user with a series of candidate positions.

\paragraph{Multi-scaling in the design space}
The distance between a local Pareto optimum and the global Pareto optima is a simple example of this.
The user can place local Pareto optima close and distant to each other, in different spatial scales.
See \cref{sec:positioning} for how we achieve this.

\paragraph{Scalability} Computing the objective functions takes $O(d v)$ steps, where $d$ is the dimension of the design space and $v$ is the number of connected components of local minima lying in the design space.

\subsection{Basin graph}
\label{sec: basin graph}

We model the effects of basins using a hypothetical solver.
Suppose we have a simple solver that improves the solution locally. 
This hypothetical solver checks the surroundings of distributed sample points and improves them by moving the sample points for an infinitesimal amount.
More precisely, given a point in the design space, the movement of our hypothetical solver is restricted so that the direction only either improves or does not change $f_i$ for any $i$.
In our model, the basin attracts this hypothetical solver, trapping them in the local Pareto set. An actual solver, of course, has sophisticated mechanisms to do a global search for the global Pareto set.
Due to this, the neighborhood of each point in the domain space is partitioned into parts, each of them can either be entered or not by this solver and leads to different connected components of the local Pareto set.
The simplistic solver continues to move its sample point, to eventually find a local (or global) Pareto set. The potential paths of the solver can be summarized into a graph that describes the different paths from a starting point. In fact, we can merge the graphs for different starting points into one, which is called the basin graph. 

s


\subsection{Design details} \label{sec:designdetails}

In a sense, the 3BC aligns with the Exploratory Landscape Analysis (ELA) \cite{mersmann2011exploratory} in the sense that the landscape will be known.
Where the ELA analyzes the landscape through a sample of points, the 3BC has a pre-defined landscape specified by the user as the basin graph we will explain later.

Our approach to realize a function from a user-specified topology has its root in visualization \cite{Weber:2007}, although they targeted scalar functions while our method targets a mapping with a 2 dimensional range.

Our approach to the benchmark problem is similar to combining spheric functions, which is routinely done in the literature.
Instead of the traditional sphere, i.e.~$\{x \in X \mid \| x \|_2 < r\}$ for some radius $r > 0$ and some spatial domain $X$, we choose a version with the Manhattan distance, $\{x \in X \mid \|x\| < r\}$.
By choosing this, we are able to specify the intersection of the spheres in $X$ more easily, which in turn lets us fine-tune the benchmark problem.

We design our benchmark problem by splitting them into a few steps:
\begin{enumerate}
\item Generate a primitive uni-objective problem
\item Generate a primitive multi-objective problem
\end{enumerate}

While each step is explained formally in later sections, here we provide some design decisions, mentioning how we fulfill the characters listed in \cref{sec:characteristicsoverview}.

As in the list above, our first step is to construct a simple uni-objective optimization problem, then proceed to a multi-objective one (which is still simplistic) and finally morphing it into a more complex, final multi-objective optimization problem, which will be fed to the solvers to be benchmarked.

\begin{figure}[tbh]
    \centering
    \includegraphics[width=0.5\columnwidth]{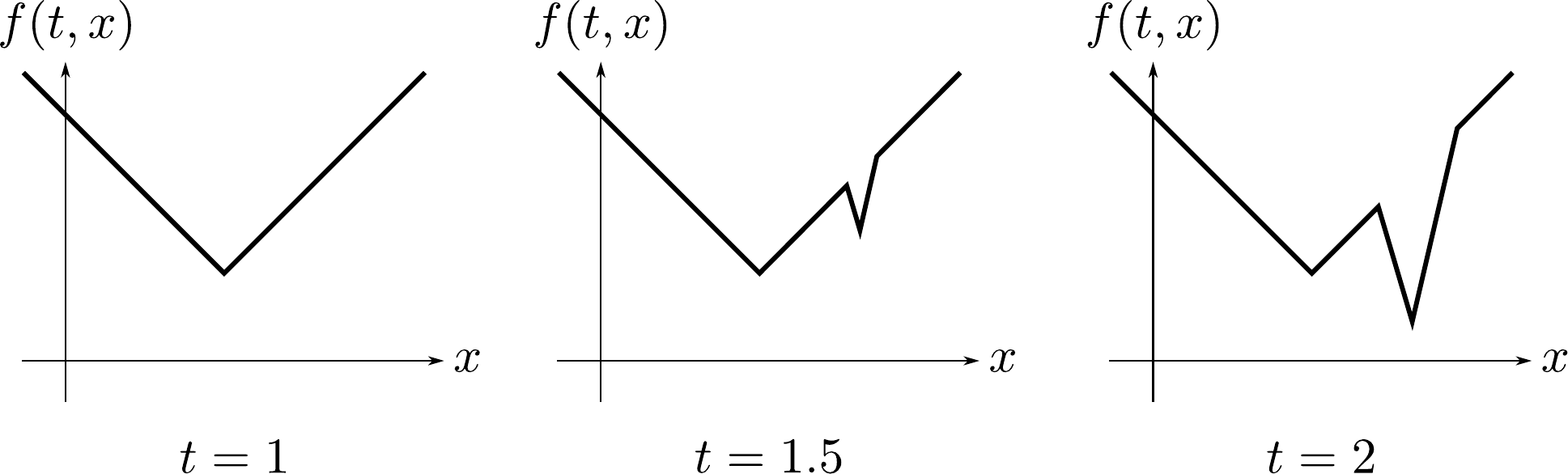}
    \caption{Graphs of $f_t$'s ($\dim X = 1$)}
    \label{fig:graphs of slice f}
\end{figure}

In step 1, we construct a uni-objective function $f$, which we call a \textit{primitive function}.
It is parameterized by a variable $t \geq 0 \in \mathbb{R}$.
Hence, $f\colon [0, \infty) \times X \rightarrow \mathbb{R}$, where $X$ is an $n$-dimensional domain specified later in \cref{eq: dfn of X} as a disk. The objective function $f$ can be regarded as a family $\{f_t\}_{t\ge 0}$ of objective functions $X\to \mathbb R$ by defining $f_t(x) \coloneqq f(t,x)$. 
At $t = 0$, the function $f_t\colon X\to\mathbb R$ is the zero function and at $t=1$, $f_t$ has a single local optimal (in fact, the minimum) solution.
As we proceed in $t$, occasionally new local minima appear (see \cref{fig:graphs of slice f}).
The user can give an ordering to the values of these local minima to decide the global minima.

\begin{figure}[tbh]
    \centering
    \includegraphics[width=0.3\columnwidth]{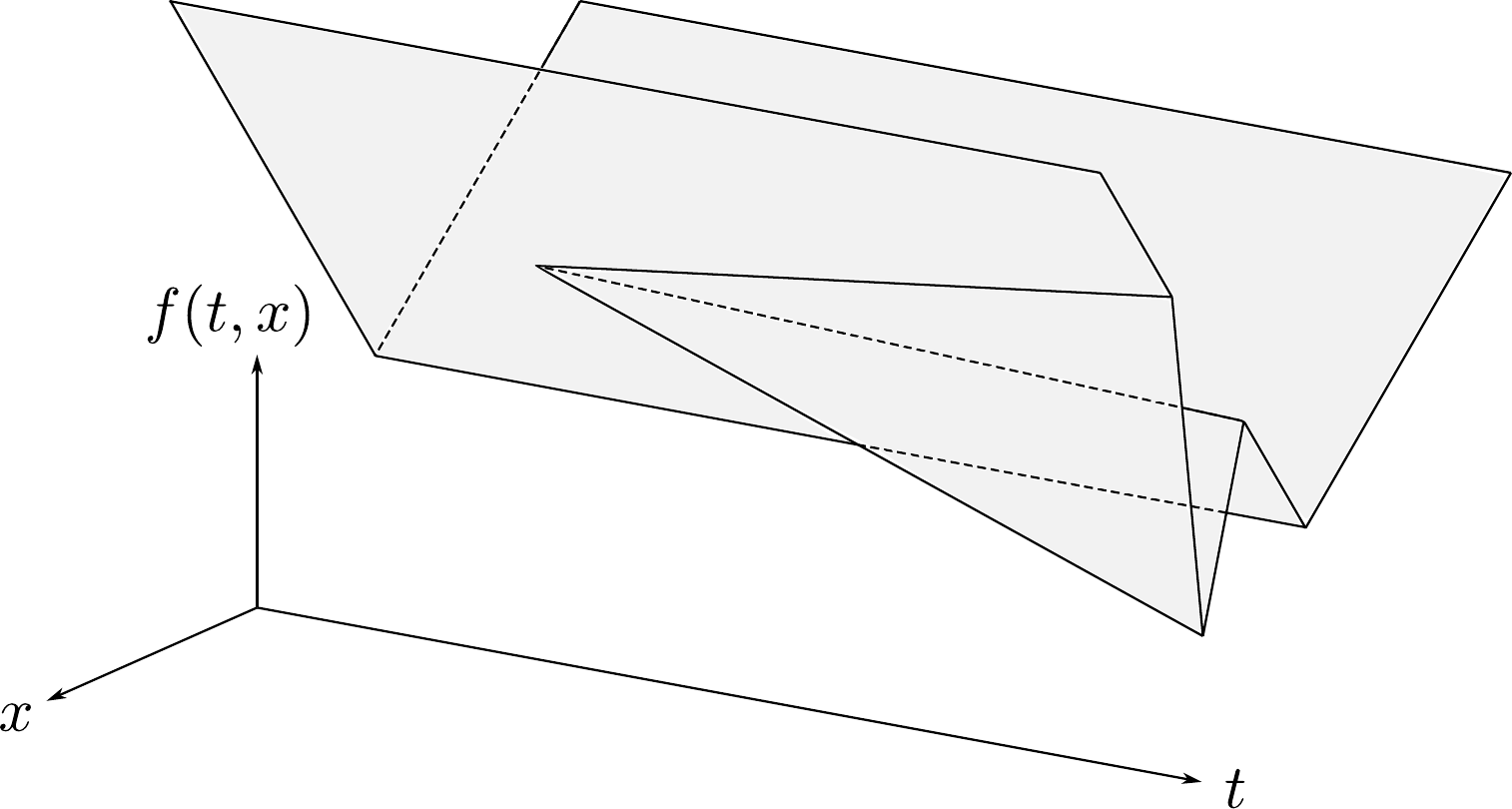}
    \caption{Graph of $f$ ($\dim X=1$)}
    \label{fig: graph of f}
\end{figure}

\begin{figure}[tbh]
    \centering
    \begin{subfigure}[b]{0.30\columnwidth}
        \centering
        \includegraphics[width=\columnwidth]{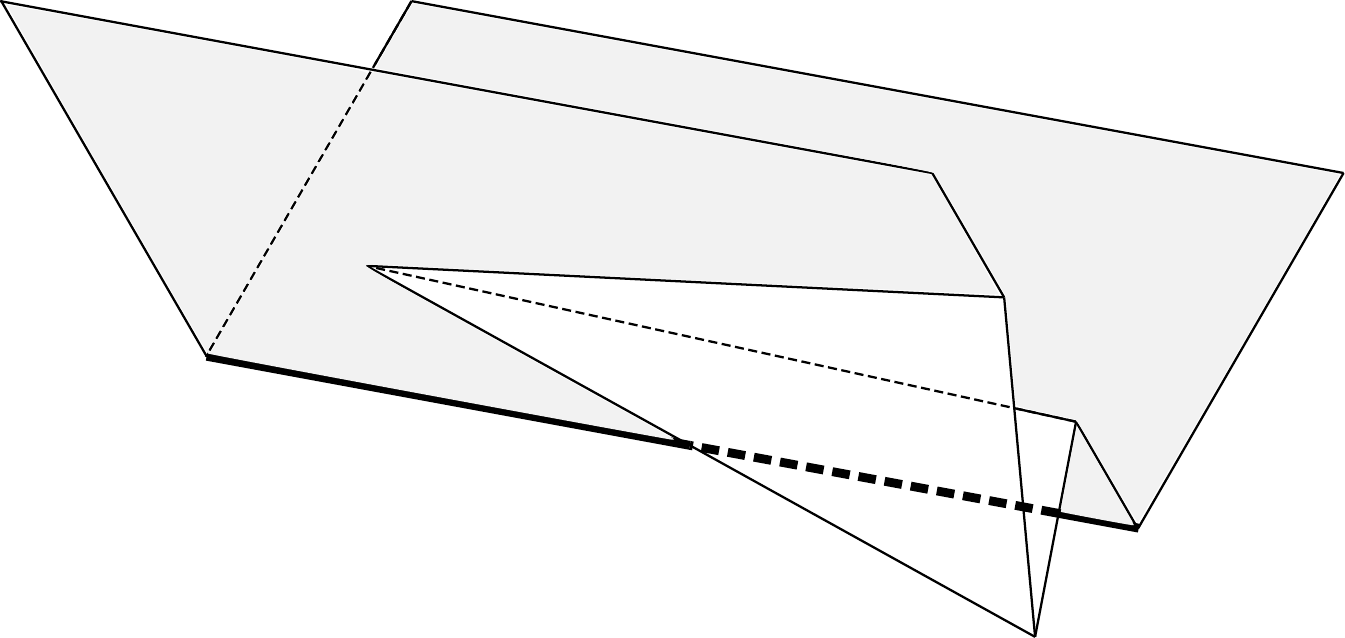}
        \caption{Basin of first minimum}
    \end{subfigure}
    \begin{subfigure}[b]{0.30\columnwidth}
        \centering
        \includegraphics[width=\columnwidth]{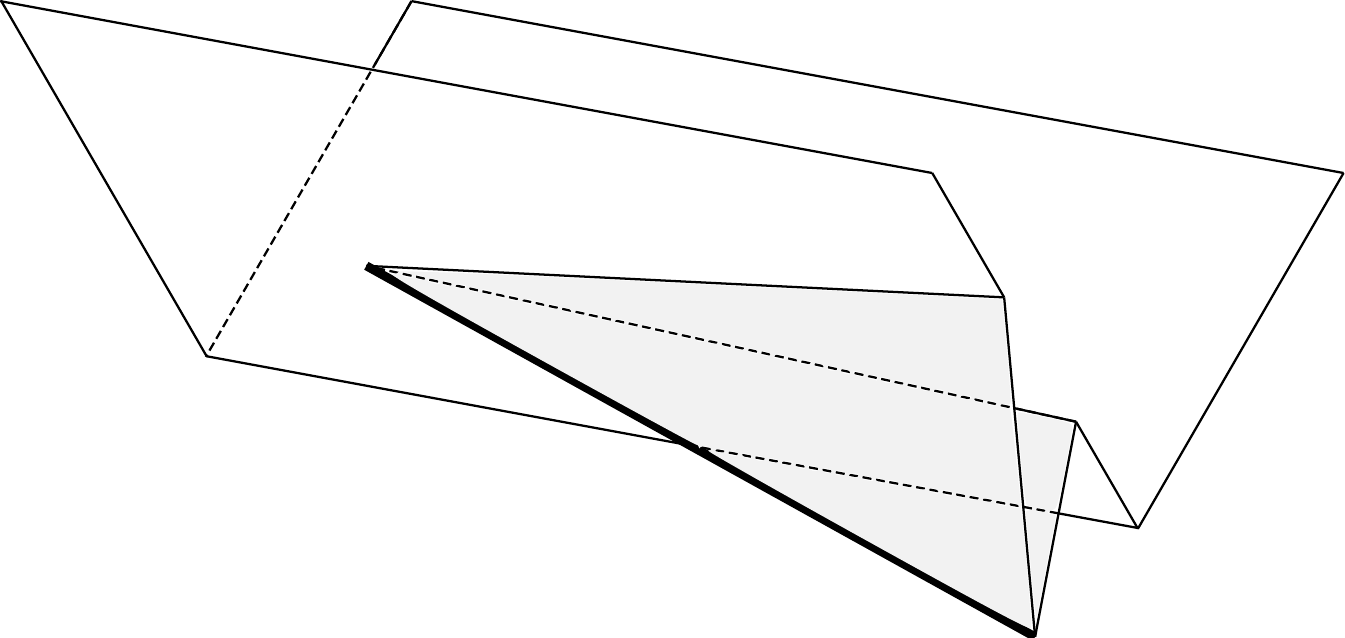}
        \caption{Basin of second minimum}
    \end{subfigure}
    \caption{Basins of local minima}
    \label{fig: basins of local minima}
\end{figure}

In the domain space of $f$, $[0, \infty) \times X$, we thus have local minima each along a line parallel to the $t$-axis (see \cref{fig: graph of f}).
These locally minimal points will be the local Pareto optima of the resulting 2-objective problem in step 2. They will be moved in step 3, although the topology of the lines of these points stays the same.
The function $f$ in \cref{fig:graphs of slice f} has two components of local minima and then their basins are described as in \cref{fig: basins of local minima}.

\begin{figure}[tbh]
    \centering
    \includegraphics[width=0.3\columnwidth]{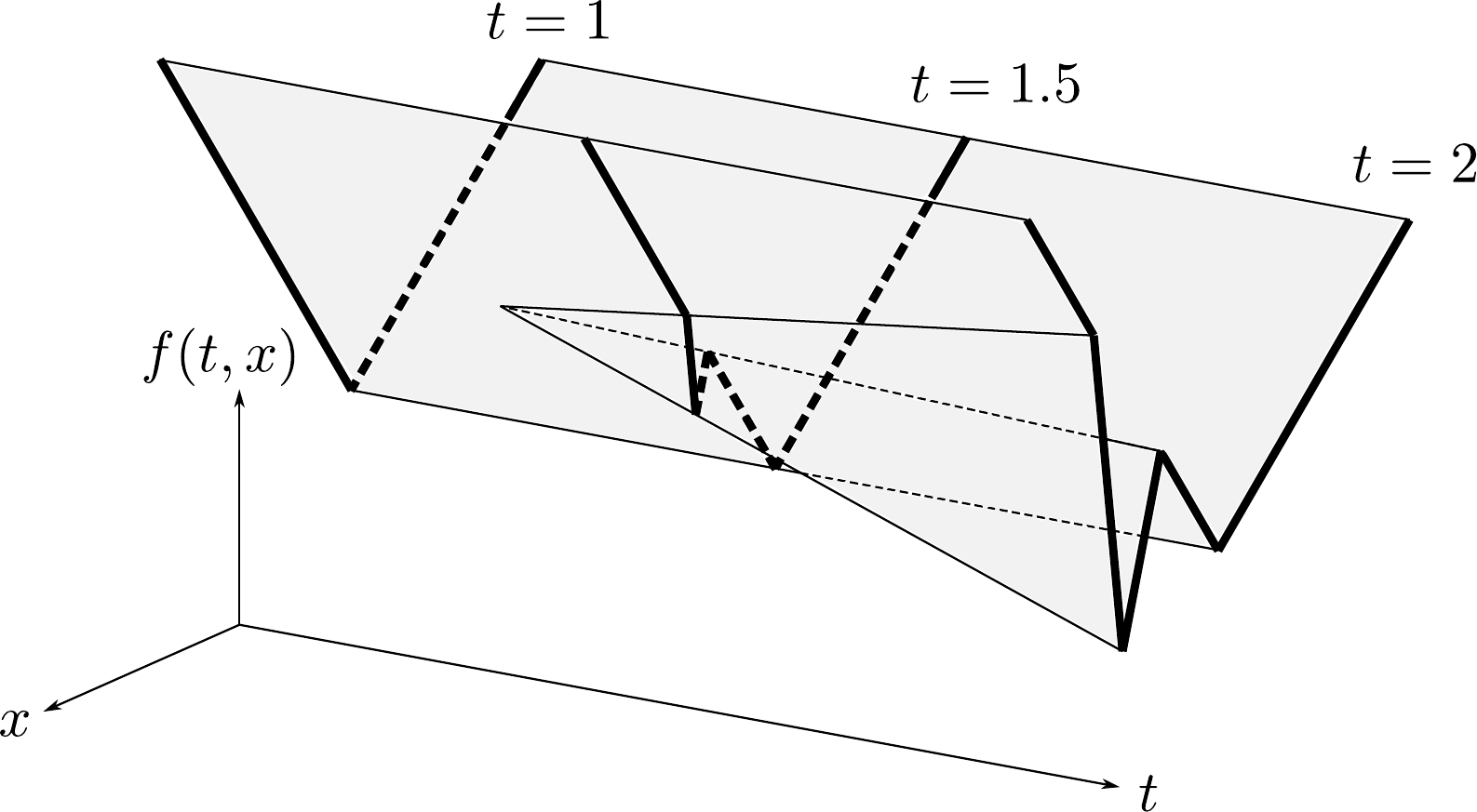}
    \caption{Family $\{f_t\}_t$ as slices of $f$ ($\dim X = 1$)}
    \label{fig: f as the family of f_t's}
\end{figure}

In step 2, we construct a 2-objective problem $[0, \infty)\times X \ni (t, x) \mapsto (t, f(t, x)) \in \mathbb R^2$. 
Then we get each $f_t$ as a slice of this 2-objective function at $t\in [0,\infty)$ (see \cref{fig: f as the family of f_t's}). 
We obtain $\tilde f$ by rotating the 2-objective function for $45^\circ$ clockwise to have the desirable Pareto set and Pareto front: 
\begin{equation}\label{eq: def of rotated objective function}
    \tilde f(t, x) \coloneqq 
    \frac{1}{\sqrt 2}
    \begin{pmatrix}
    1 & 1 \\ -1 & 1
    \end{pmatrix} 
    \begin{pmatrix}
    t \\ f(t, x)
    \end{pmatrix}.
\end{equation}

The final problem is obtained as follows. 
We can make the Pareto set of the problem in step 2 more complicated, which would reduce the performance of solves. 
We can also deform the shape of the Pareto front, which has an effect on the performance of algorithms because a different algorithm is good at a different shape of the Pareto front.

Characteristics of our benchmark problem can be explained by the condition on the critical values $m_s$'s of $f$, 
which will appear in \cref{sec: appendix dealing with separability}. 
For example, Separability is achieved if $m_\emptyset$ is the minimum of $m_s$'s. 
The critical values $m_s$'s also determine the shape of Pareto front and especially the Pareto front is completely convex when $m_\emptyset$ is the minimum of $m_s$'s.


\section{Overview of the primitive function $f$}

\subsection{Property of the primitive function $f$}

\begin{figure}[bthp]
    \centering
    \includegraphics[width=0.5\columnwidth]{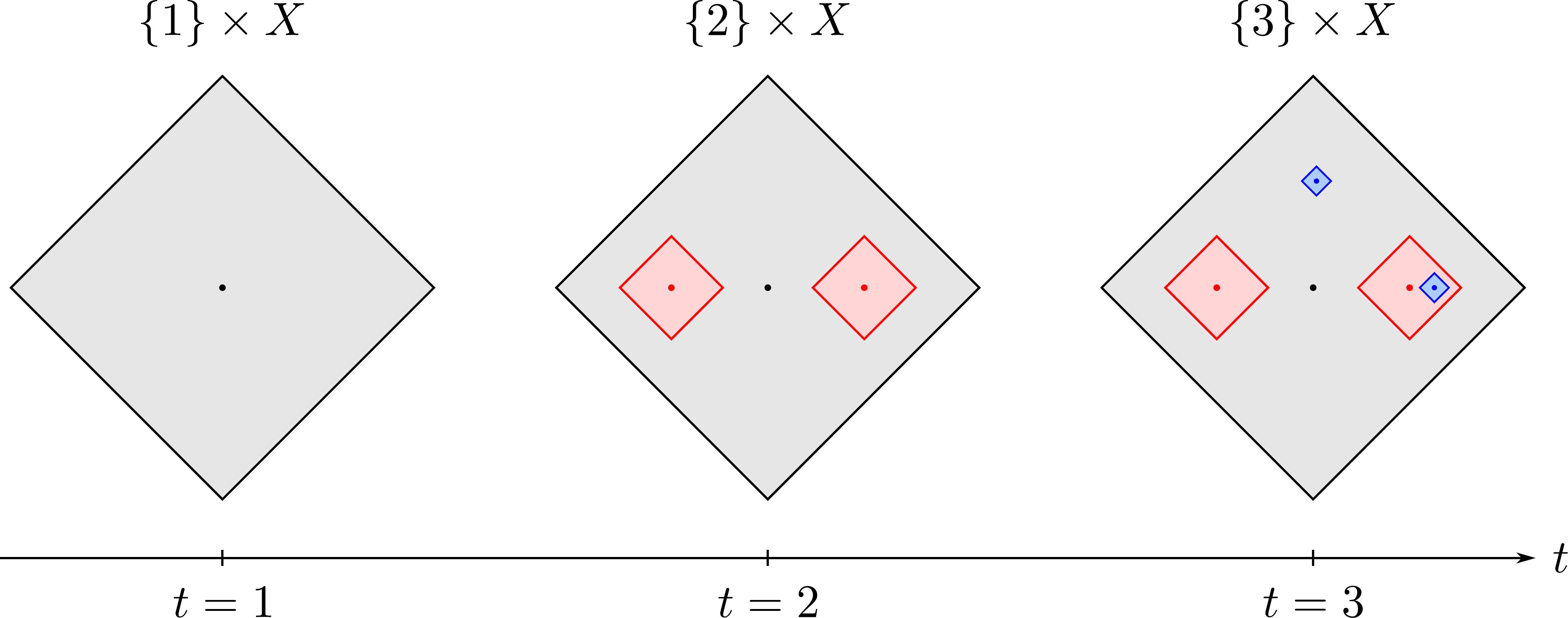}
    \caption{An example of diamonds in $X$ ($n=\dim X=2$)}
    \label{fig: diamonds as t increases}
\end{figure}

As mentioned in \cref{sec:designdetails}, we regard $f\colon [0, \infty)\times X\to \mathbb R$ as a family $\{f_t\}_{t\ge 0}$ of single-objective functions on $X$ such that local minima are added as $t$ increases. 
In this section, we explain what properties we request for $f$. 
First, the domain $X$ is defined to be the unit ball in $\mathbb R^n$ whose shape is like a diamond, that is, 
\begin{equation}\label{eq: dfn of X}
    X = \{x\in \mathbb R^n\mid \|x\|\le 1\}, 
\end{equation}
where $\|x\| = \sum_{i=1}^n|x_i|$ is the norm of $x$. 
Then $X$ is the black diamond on the left of \cref{fig: diamonds as t increases}. 
When $t=1$, we construct $f$ so that it has only one local minimum at the origin. 
When $t$ increases and reaches 2, the function $f_t$ can have more local minima other than the origin. 
Each local minimum has its diamond in which the function has only one local minimum, for example, the two red diamonds in the middle figure in \cref{fig: diamonds as t increases}). The specific definition of the diamond is given in \cref{sec:positioning}, in particular in \eqref{eq: def of the diamond}.  
Similarly, $f_3$ may have more local minima other than those of $f_2$.

\subsection{Positioning local minima and their diamond} \label{sec:positioning}

\begin{figure}[bth]
    \centering
    \includegraphics[width=0.35\columnwidth]{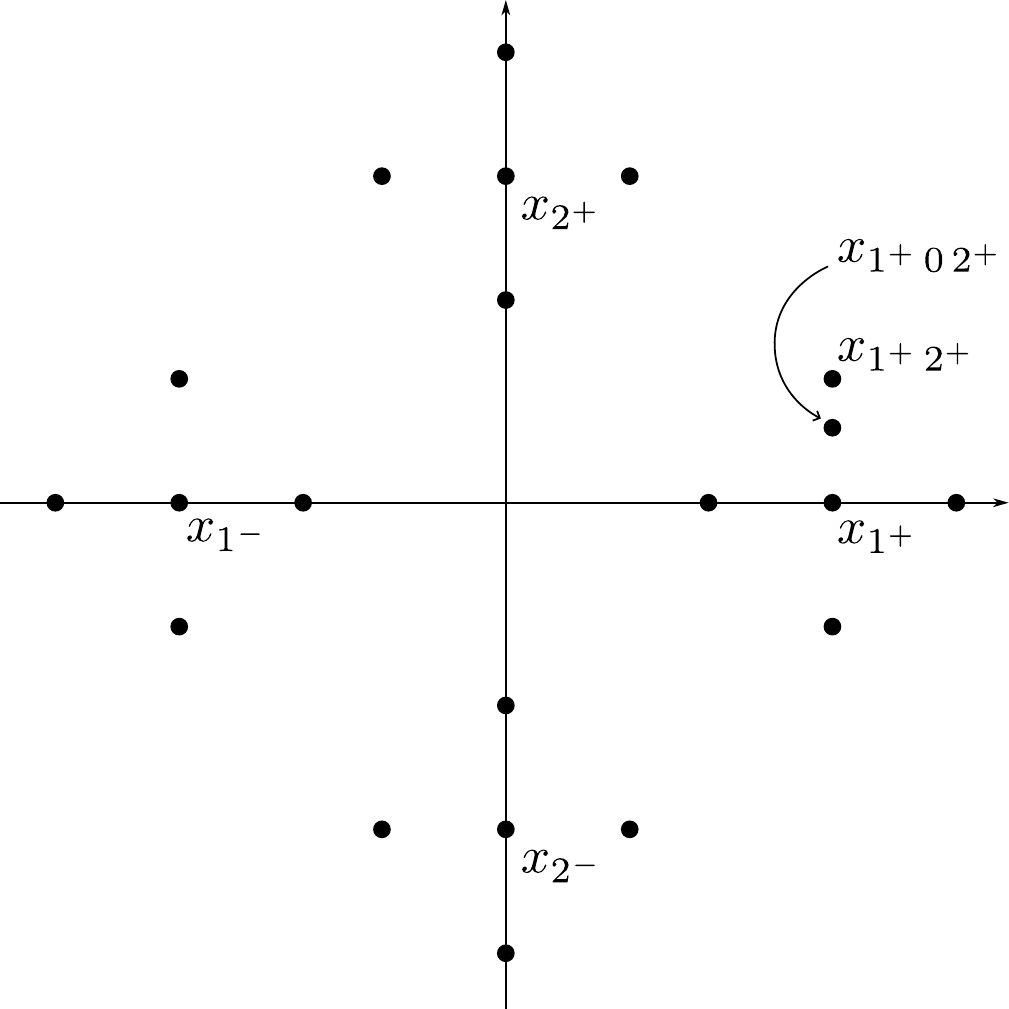}
\caption{Position of local minima ($\operatorname{dim}X = 2$)}
\label{fig:sequenceexample}
\end{figure}

Local minima of the primitive function $f$ can be chosen from predefined candidate positions.
To allow the benchmark to model multiple scales in the design space, we place the position in a fractal structure (see \cref{fig:sequenceexample}).
The position of each candidate point is coded as a sequence $\movements$ of \emph{movements}.
We start with the origin of the coordinate system and follow the movements coded in $\movements$ to specify the candidate position $\candidate{\movements}$.
Each $j$-th movement can be taken along an axis in $X$, and it can either go backwards or forwards for length
$2/4^j$. 
In addition, a movement can also stay at the position.
We denote each movement along $i$-th axis backwards or forwards as $i^-$ and $i^+$, respectively. In addition, a ``movement'' can also be a \emph{stay}, indicated by a 0,
and so we can write down a sequence $\movements$ of movements as, e.g.,
\begin{equation} \label{eq:movements}
	 \movements \coloneqq 1^{+}\,0\,2^+.
\end{equation}
This example results at the position illustrated in \cref{fig:sequenceexample}.
We require the last movement of any sequence not to be $0$. 
In the example indicated in \cref{fig: diamonds as t increases}, the appearing movements are $s=\emptyset$, $1^+$, $1^-$, $1^{+}\,1^+$, and $0\,2^+$. 

We denote the number of movements in $\movements$ as $|\movements|$ ($= 3$ for \cref{eq:movements}).
We allow $\movements = \emptyset$, which denotes a local minimum at the origin.
Indeed, this local minimum is a mandatory user input in our benchmark for the sake of simpler mathematics.
The user specifies the position of the local minimum with a sequence of movements.

For each $s$, its diamond $\dia{s}\subset X$ is defined by 
\begin{equation}\label{eq: def of the diamond}
    \dia{s} = \{x\in X\mid \|x-x_s\|\le 1/4^{|s|}\}. 
\end{equation}
Then diamonds of movements of the same length are all disjoint. 
Let $s$ be a movement and consider its ``child" $s' = sj^{\pm}$ for some $1\le j\le n$. 
Since $\|x_s-x_{s'}\| = 2/4^{|s|+1}$, the diamond of 
$s'$ is a subset of that of $s$, which means the diamonds are nested.

\subsection{Parameter-family of local optima as local Pareto set} \label{sec:localoptima}

\begin{figure}[tbh]
    \centering
    \includegraphics[width=0.6\textwidth]{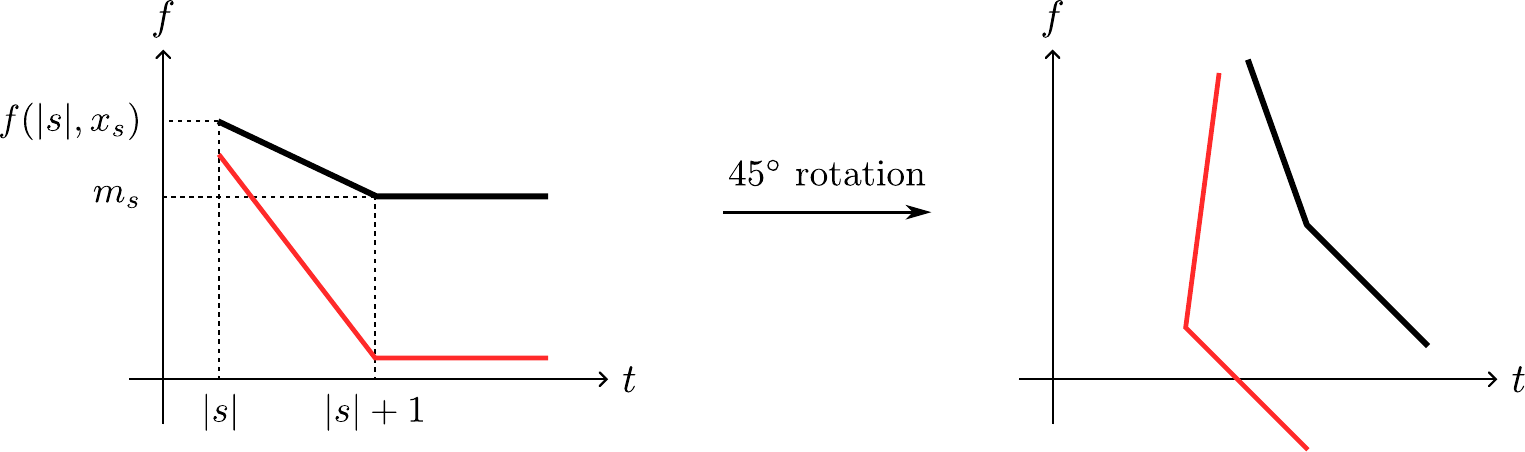}
    \caption{Local minima of $f$ and a local Pareto front of $\tilde f$}
    \label{fig: condition of Pareto front}
\end{figure}

By increasing $t$, we give rise to new local minima of $f$, which will be the local Pareto optima of the 3BC.
We assume $f$ to be continuous for the sake of the simplicity of the mathematics.
Each local minimum at position $\candidate{\movements}$ emerges right after $t = |\movements|$, having the local minimum value $M_\movements(t)$ and reach the user-specified value $m_\movements$ at $t = |\movements| + 1$.
After this $t$, the local minimum stays constant at this value, and will stay being a local minimum. 
The precise definition of $M_s(t)$ appears in \cref{subsec: For points inside a diamond} and its graph is on the left side of \cref{fig: condition of Pareto front}. 
As in \cref{fig: condition of Pareto front}, the set $\{(t, x_s)\mid t\ge |s|\}$, which is the component of local minima of $f$ for node $s$, is also a connected component of the local Pareto sets of $\tilde f$ if 
\[
    m_\movements - f(|\movements|, x_\movements) > -1. 
\]
Otherwise the connected component is $\{(t, x_s)\mid t\ge |s|+1\}$, that is, each point in the set $\{(t, x_s) \mid |s|\le t\le |s|+1\}$ is not a local optimum (see red lines in \cref{fig: condition of Pareto front}). Note that $m_s$ must satisfy more preliminary condition~\eqref{eq: condition on m_s} appearing later to be a local minimum of $f$. 





\subsection{Specifying basins using the basin graph} \label{sec:inputtree}


Recall that the user specifies the position of a local minimum using a sequence of movements as in \cref{eq:movements}.
This defines a tree structure, which we call the \textit{basin graph}. 
Our basin graph is similar to the \emph{reachability graph} \cite{Bormann2020}
Each node corresponds to a movement, which in turn corresponds to a connected component of the local minima of $f$. – the root node of this tree corresponds to the local minimum without a movement, i.e.~$\movements = \emptyset$.
Let $\movements$ be the sequence of movements for a node in the tree.
This node is connected to a child node with a link if and only if
    (i) $\movements$ matches the first $|\movements|$ items in the sequence of the child node, and
    (ii) the child node is not an offspring of another child of this parent node. 
We call the local minimum corresponding to the parent and child basins simply a \emph{parent} or a \emph{child}.

By using the basin graph, the basin of each node can be defined as follows. 
For a node $s$, from its diamond $\dia{s}$ we obtain
\[
    \bigcup_{t\ge |s|} \{t\}\times \dia{s}, 
\]
which shapes a cone since the radius of the diamond increases monotonically. 
Then, the basin of a node is defined as its cone subtracted by cones of all children of the node. 
\begin{figure}[tbh]
    \centering
    \includegraphics[width=0.2\columnwidth]{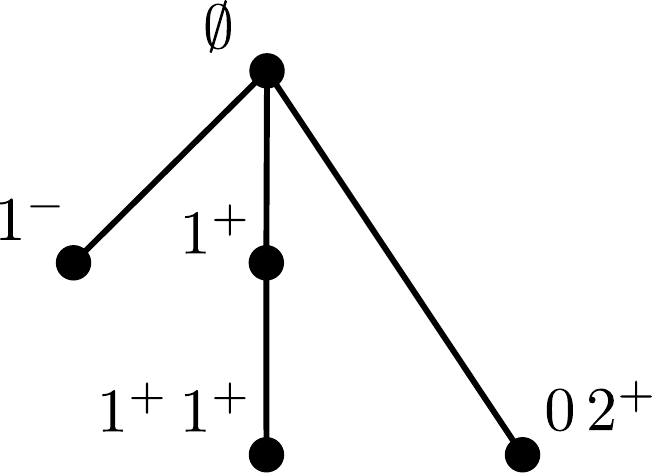}
    \caption{The basin graph corresponding \cref{fig: diamonds as t increases}}
    \label{fig: an example of the basin graph}
\end{figure}



\section{Defining and computing the primitive function $f$} \label{sec:compute}


We define the objective function $f$ for $t = 0$ and then $t > |\movements|$ inside the diamond
\begin{equation*}
    \dia{s} = \{x\in X\mid \|x-x_s\|\le 1/4^{|s|}\} \tag{\ref{eq: def of the diamond}}
\end{equation*}
of each $\movements$.
As the domain of $f$ is restricted to the largest diamond (\cref{eq: dfn of X}) for at any $t$, $f$ is indeed assigned value at any point.
If a point is inside two diamonds of sequences $\movements_1$ and $\movements_2$, we choose $f$ defined by the cone of the younger $s_i$ (i.e.~the one with larger $|\movements_i|$).

In the next two sections, we will define $f \colon [0, \infty)\times X\to \mathbb R$ by inductively constructing $f \colon (\tau, \tau+1]\times X\to \mathbb R$ for $\tau=0,1, 2, \dots$

\subsection{At $t = 0$}

{\color{red}
}

Let $x$ be a point in $X$.
We define $f(t, x)$ as
\[
\begin{aligned}
f(t, x) \coloneqq 0 && \text{for } t = 0.
\end{aligned}
\]

\subsection{At $\tau < t \leq \tau + 1$.} \label{subsec: For points inside a diamond}

Assuming that $f$ is defined for $\{\tau\}\times X$, where $\tau$ is a non-negative integer, we construct $f$ for $(t, x) \in (\tau, \tau+1]\times X$ in the following manner.
We consider the sequences of movements whose length is $\tau$ (if they exist in the user input). Let $s$ (with $|s| = \tau$) be such one. If there is none, $f(t,x) \coloneqq f(\tau,x)$ for any $x \in X$.

In the remainder, we consider the situation when such a sequence $s$ exists.
Let $S_\tau$ be the set of such sequences. If so, $f(t, x)$ for $(t, x) \in (\tau, \tau+1]\times X$ is defined to be

\begin{equation} \label{eq:defF}
    f(t, x) 
    \coloneqq
    \min_{|s| < t} \{ g_s(t, x) \}.
\end{equation}
(see 
\cref{fig:graphs of slice f}
). 
The function $g_s$ is explained in the rest of this subsection.

\begin{figure}[t]
    \centering
    \includegraphics[width=0.3\textwidth]{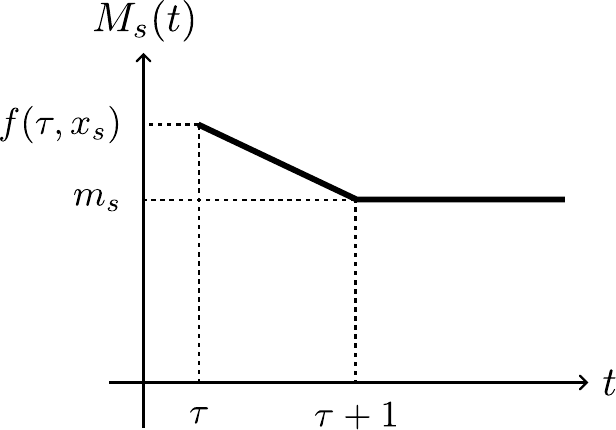}
    \caption{Value of $f$ at a local minimum $x_s$.}
    \label{fig:fdshfdsf}
\end{figure}

As shown in \cref{fig:fdshfdsf}, the locally minimum value $M_\movements(t)$ at position $\candidate{\movements}$ firstly appears right after $t = \tau$, decreasing its value to the user-given $\con_s$ via linear interpolation, till $t = \tau + 1$:
\begin{equation*}
    M_s(t) \coloneqq (1-\DeltaT)\ f(\tau, \cp_s) + \DeltaT \ \con_s,
\end{equation*}
where $\DeltaT \coloneqq \min(t-\tau,1)$ and $t \geq \tau$.
$\min$ is taken to stop the growth of $\func_s(t,x)$ at $t = \tau+1$ and we impose the following condition on $m_s$'s specified by the user. 
\begin{equation}\label{eq: condition on m_s}
    m_\movements < f(|\movements|, x_\movements).
\end{equation}

Function value $\func_s(t, x)$ at parameter $t$ on a point $x$ around the local minimum point $\candidate{\movements}$ is given by climbing upwards from $\candidate{\movements}$, following the gradient vector $\nabla \func_s(x)$:

\begin{equation} \label{eq:func_s}
    \func_s(t, x) 
    \coloneqq M_s(t) + \nabla \func_s(x) \cdot \Delta x,
\end{equation}
where $\Delta x \coloneqq x - \candidate{\movements}$ is a vector expressing the deviation from the locally minimal point.

\begin{figure}[t]
    \centering
    \includegraphics[width=0.5\textwidth]{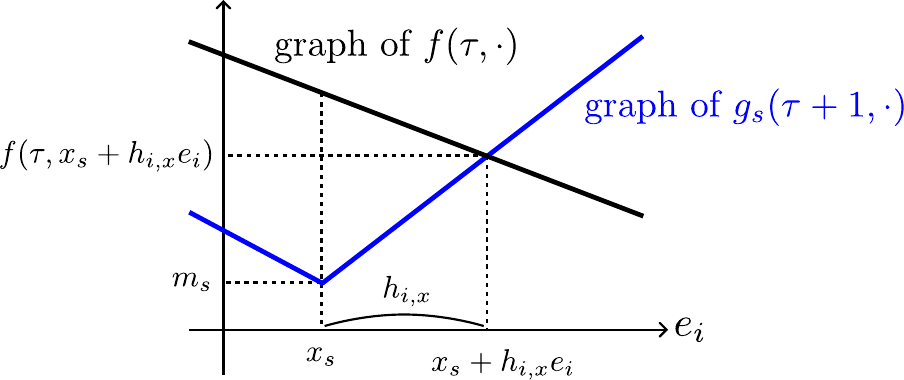}
    \caption{Definition of $\nabla g_s(x)$}
    \label{fig: gradient of g_s}
\end{figure}

To decide $\nabla \func_s(x)$, we compute the intersection points between the graphs of $\movements$ and the parent of $s$ (see \cref{fig: gradient of g_s}).
Note that $\nabla \func_s(x)$ is independent of $t$.

We require any, say $i$-th, component of $\nabla \func_s(x)$ to match
\begin{equation} \label{eq:nabla_g_future}
	\frac {
				f(\tau + 1,\ \candidate{\movements} + \ \movement{i}{x}\ e_i)
				-
				f(\tau + 1,\ \candidate{\movements})
			}{
				\movement{i}{x}
			} \eqqcolon \nabla \func_s^i(x),
\end{equation}
which follows from the definition of a differential.
Here, $e_i$ is the $i$-th unit vector and  ${\candidate{\movements}} + \movement{i}{x}\ e_i$ gives location of $x \in X$ at which the graph of $g_\movements(t, x)$ intersects $f(\tau, \bullet)$:
\begin{equation*}
	\movement{i}{x}\coloneqq\sgn(i, x)\frac{1}{4^\tau} \ \in \mathbb{R}.
\end{equation*}
$\sgn(i, x)$ tells the side of the diamond on which point $x$ lies along the $i$-th axis, i.e.~the sign of $x^i - x^i_s$.
If this subtraction tends to 0 (or, $x^i \rightarrow x^i_s$), $\nabla g^i_s(x) \times (x^i - x^i_s)$ in \cref{eq:func_s} approaches 0, regardless of the value of $\sgn$. In our implementation, we thus set $\sgn(i, x_s)$ simply to 1.

In fact, regarding \cref{eq:nabla_g_future}, $f(t,\ \candidate{\movements} + \ \movement{i}{x}\ e_i)$
is constant as the $\candidate{\movements} + \ \movement{i}{x}\ e_i $ is outside the diamond of $s$. In addition, $f(\tau + 1,\ \candidate{\movements}) = m_s$ because $\candidate{\movements}$ is the local minimum point. Hence,
\begin{equation*}
    \nabla \func_s^i(x) \equiv
	\frac {
				f(\tau,\ \candidate{\movements} + \ \movement{i}{x}\ e_i)
				-
				m_s
			}{
				\movement{i}{x}
			},
\end{equation*}
and we can determine $g_s$ using \cref{eq:func_s}.

Finally, we decide $f \colon (\tau, \tau+1]\times X\to \mathbb R$:
\begin{equation*}
    f(t, x) 
    \coloneqq
    \min_{|s| < t}  g_s(t, x)  
    =
    \min\left(\{f(\tau, x)\}\cup \{g_s(t, x) \mid |s| = \tau\}\right). 
\end{equation*}
The first equality is due to \cref{eq:defF}.
\begin{align*}
    f(t,x)
    &=
    \min_{|s|\le \tau} g_s(t,x) \nonumber \\
    &=
    \min (\{g_s(t,x)\mid |s|\le \tau-1\}\cup \{g_s(t,x)\mid |s|=\tau\}) \nonumber \\
    &=
    \min (\{g_s(\tau,x)\mid |s|\le \tau-1\}\cup \{g_s(t,x)\mid |s|=\tau\}) \nonumber \\
    &=
    \min (\{f(\tau, x)\}\cup \{g_s(t,x)\mid |s|=\tau\})
\end{align*}

As shown in \cref{fig: gradient of g_s}, and by equation~\ref{eq:func_s}, for $t\in (\tau, \tau+1]$ and a sequence $s$ with $|s|=\tau$, $g_s(t, x) \le f(\tau, x)$ only in the diamond of $s$. The function $f_\tau$ does not have $x_s$ as a critical point, but $g_s$ does for every $t$, so minimization above can make $x_s$ a critical point after $t\ge \tau$. It is our desirable condition. 

\section{Basin-Wise IGDX} \label{sec: basin-wise IGDX}

It is useful to measure whether each solver finds the optimal solutions. 
We can use various performance indicators such as GD~\cite{van1999multiobjective}, IGD~\cite{coello2004study}, GDX~\cite{schutze2011computing}, and IGDX~\cite{schutze2011computing, zhou2009approximating} to do it since the Pareto set and front of our benchmark problem are known. 
Hypervolume~\cite{zitzler1998multiobjective} is also useful. 

3BC is designed to study the behavior of Evolutionary Algorithms (EAs) that change at each generation, reacting to the basins.
As a result, we study the population at each generation of different EAs such as GDE, IBEA etc.
In other words, we do not use the external archive for our evaluation of the 3BC – such a strategy is more suited for studying the overall performance of EAs instead of the reaction to the basins for the situation EAs is in at the temporal instant.

It is also useful to measure whether solutions obtained by a solver is close to the local Pareto sets or not. 
The IGDX is also useful here by taking $P^*$ (this notation appears below) as a subset of the local Pareto sets. 
However, the IGDX can take a good value even when solutions obtained by an algorithm are distributed all over the place.
It can be improved by combining the IGDX above with the GDX since the more solutions obtained by an algorithm away from the local optima, the worse the GDX. 

These indicators cannot measure how much the solver is misled by each local optimum along the way, so another performance indicator is needed to measure that. 
One conventional way to measure the extent to which a solver is trapped in a local Pareto solution is to use $x$-means. This is simply a measure of how many clusters it is appropriate to divide the solver's solution into. 
In contrast, the 3BC defines the basin of attraction for each node (connected component of the local Pareto set). 
We can be use it to define an indicator for each basin, which we call a \textit{basin-wise} indicator for each basin. 
In particular, \textit{basin-wise} IGDX for each basin of attraction can measure how much a solver's solution falls into its basin. 

Let $P$ be the solution by a solver and $P^*$ be the (global) Pareto set. 
Then original IGDX is defined by 
\[
    \operatorname{IGDX}(P^*, P) = \dfrac{\Sigma_{v\in P^*}\,d(v, P)}{|P^*|}, 
\]
where $d(v, P)$ is the Euclidian distance in the decision space. 
Replacing $P^*$ by a local Pareto optimum may define a IGDX for the basin but this indicator has an undesirable property: 
The situation where the solutions $P$ is stuck in other basins has an unintended effect on the value of the indicator. 
This is due to the fact that the entire solution is taken as $P$, and it would be a good indicator if we focus only on the solutions related to the basin. In general, it is not possible to do so, but 3BC has a basin of attraction, which can be defined. 
That is, we define basin-wise IGDX for a node $P^*$ by the value
\[
    \operatorname{IGDX}(P^*, B) = \dfrac{\Sigma_{v\in P^*}\,d(v, B)}{|P^*|}, 
\]
where $B$ is the basin of $P^*$.

\section{Evaluation} \label{sec: evaluation}

\begin{figure}[tbh]
    \centering
    \begin{subfigure}[b]{0.45\columnwidth}
        \centering
        \includegraphics[width=0.65\columnwidth]{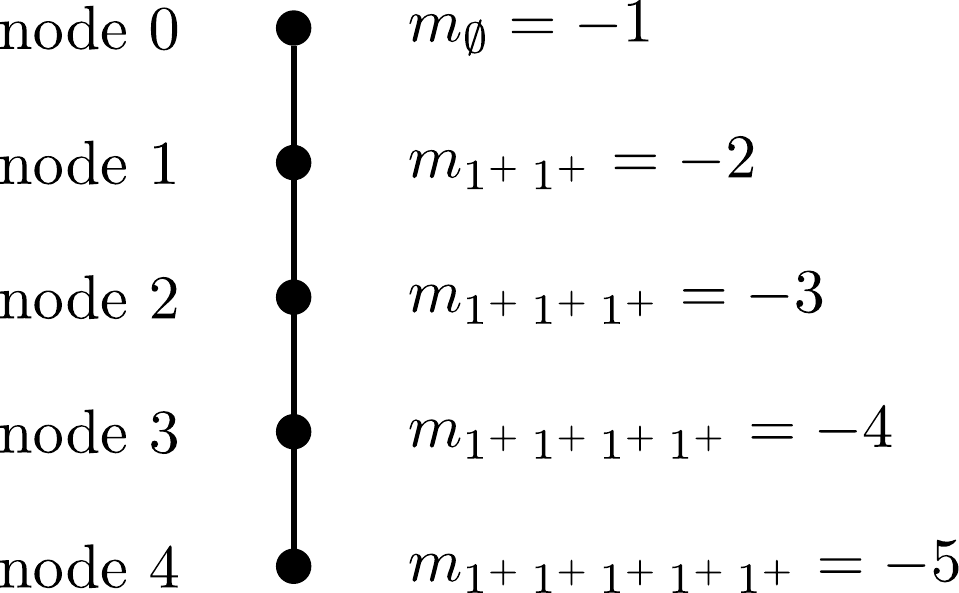}
        \caption{A depth-base basin graph}
        \label{fig: basin graph of depth 2}
    \end{subfigure}
    \begin{subfigure}[b]{0.45\columnwidth}
        \centering
        \includegraphics[width=0.9\columnwidth]{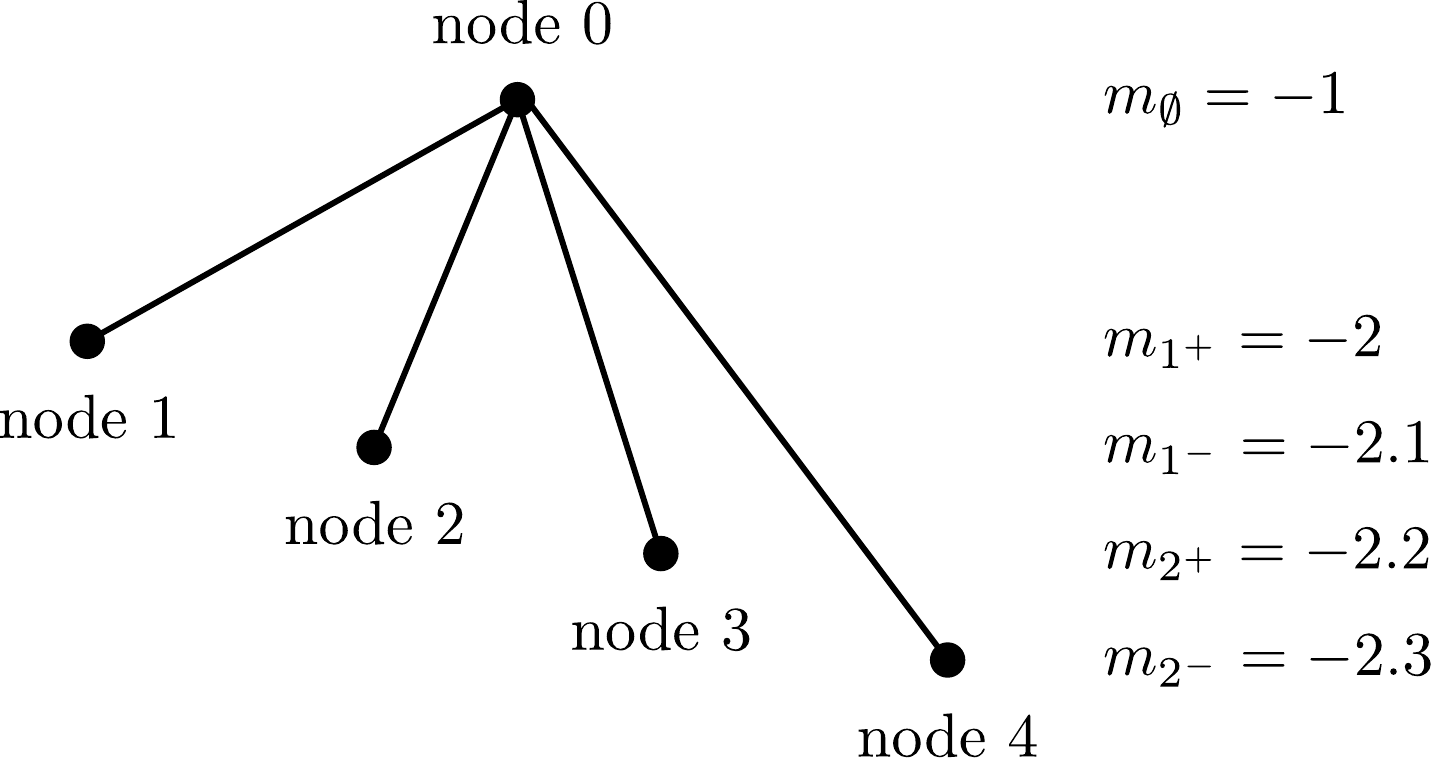}
        \caption{A breadth-base basin graph}
        \label{fig: basin graph of breadth 1}
    \end{subfigure}
    \caption{Examples of the basin graph}
    \label{fig: examples of the basin graph}
\end{figure}

When an algorithm solves the 3BC, the values of basin-wise IGDX in each generation construct a line graph. 
For example, consider a 3BC having a basin graph in \cref{fig: basin graph of depth 2}, which we call a depth-base basin graph.
This 3BC has five nodes $\emptyset$, $1^+\,1^+$, $1^+\,1^+\,1^+$, $1^+\,1^+\,1^+\,1^+$, and $1^+\,1^+\,1^+\,1^+\,1^+$. 

\begin{figure}[tbh]
    \centering
    \begin{subfigure}[b]{0.18\columnwidth}
        \centering
        \includegraphics[width = \columnwidth]{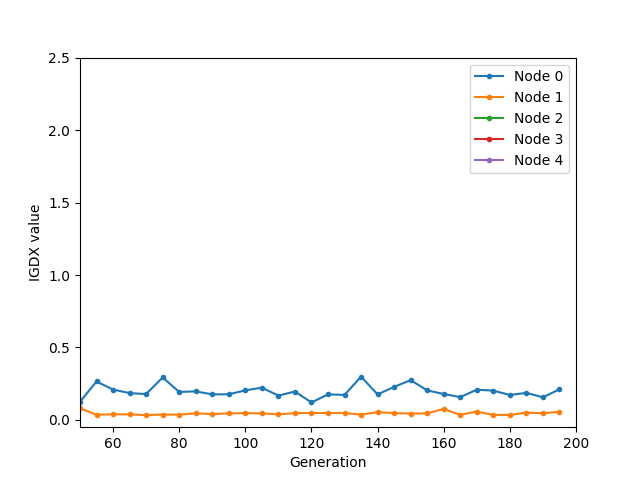}
        \caption{GDE}
    \end{subfigure}
    \begin{subfigure}[b]{0.18\columnwidth}
        \centering
        \includegraphics[width = \columnwidth]{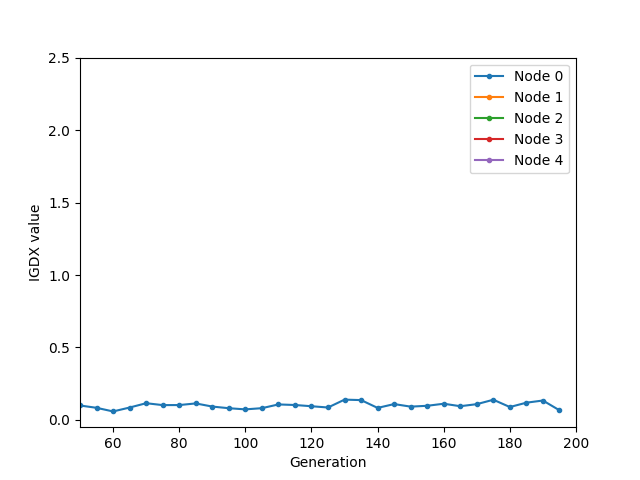}
        \caption{IBEA}
    \end{subfigure}
    \begin{subfigure}[b]{0.18\columnwidth}
        \centering
        \includegraphics[width = \columnwidth]{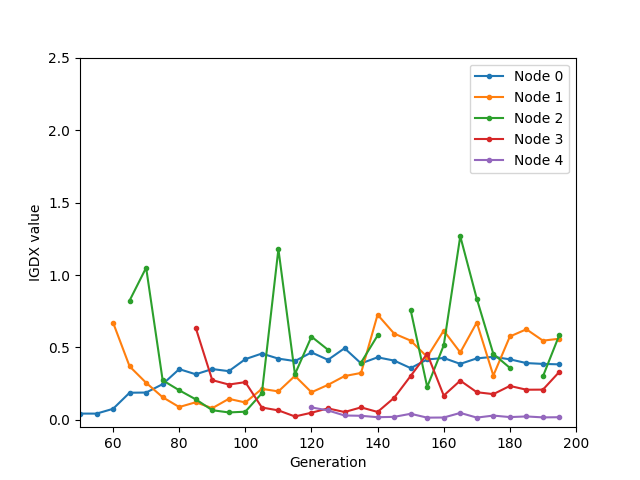}
        \caption{MOEA/D}
    \end{subfigure}
    \begin{subfigure}[b]{0.18\columnwidth}
        \centering
        \includegraphics[width = \columnwidth]{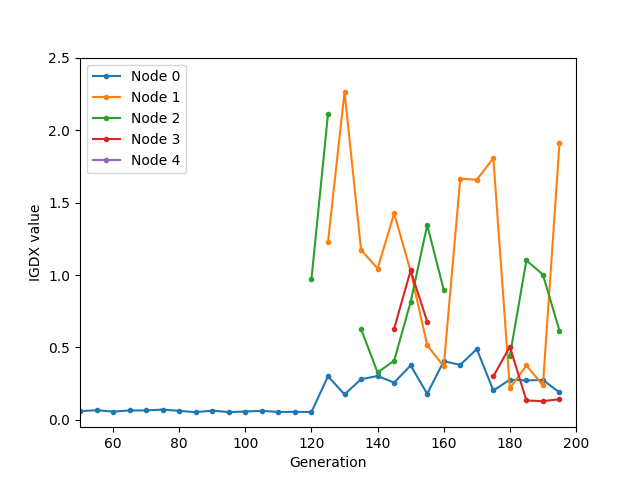}
        \caption{NSGA-II}
    \end{subfigure}
    \begin{subfigure}[b]{0.18\columnwidth}
        \centering
        \includegraphics[width = \columnwidth]{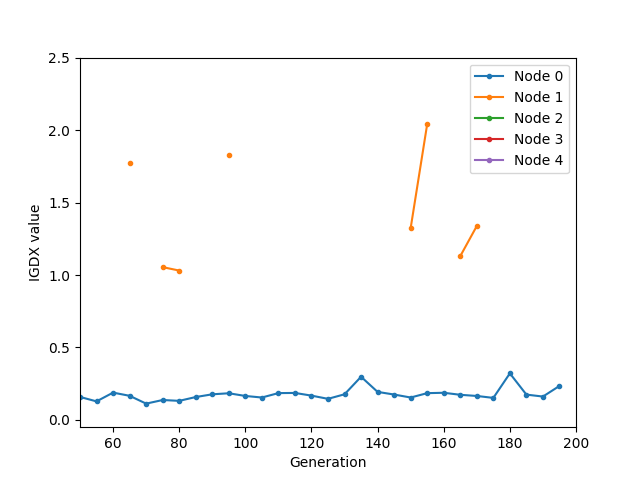}
        \caption{OMOPSO}
    \end{subfigure}
    \caption{Basin-wise IGDX for each basin of attraction (for a depth-base basin graph)}
    \label{fig:basin-wiseIgdx}
\end{figure}

Let the values $m_s$'s be $-1$, $-2$, $-3$, $-4$, and $-5$. 
The line graphs of the basin-wise IGDX are as shown in \cref{fig:basin-wiseIgdx}, where the solvers are GDE, IBEA, MOEA/D, NSGA-II, and OMOPSO.

\begin{figure}[tbh]
    \centering
    \begin{subfigure}[b]{0.18\columnwidth}
        \centering
        \includegraphics[width = \columnwidth]{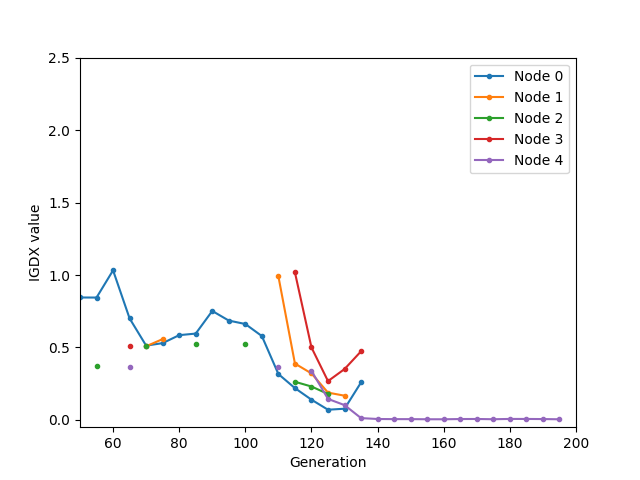}
        \caption{GDE}
    \end{subfigure}
    \begin{subfigure}[b]{0.18\columnwidth}
        \centering
        \includegraphics[width = \columnwidth]{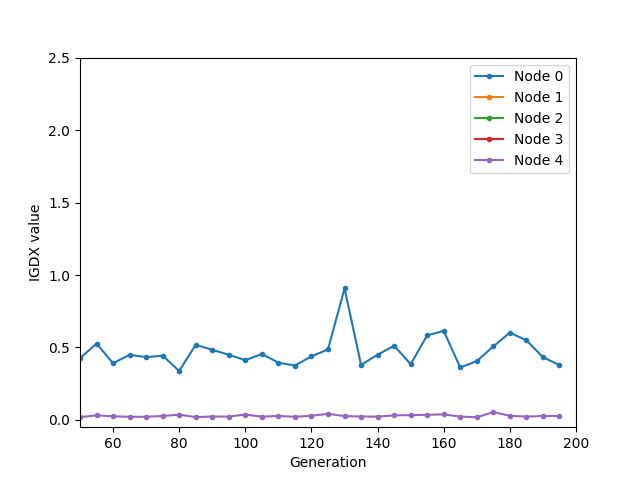}
        \caption{IBEA}
    \end{subfigure}
    \begin{subfigure}[b]{0.18\columnwidth}
        \centering
        \includegraphics[width = \columnwidth]{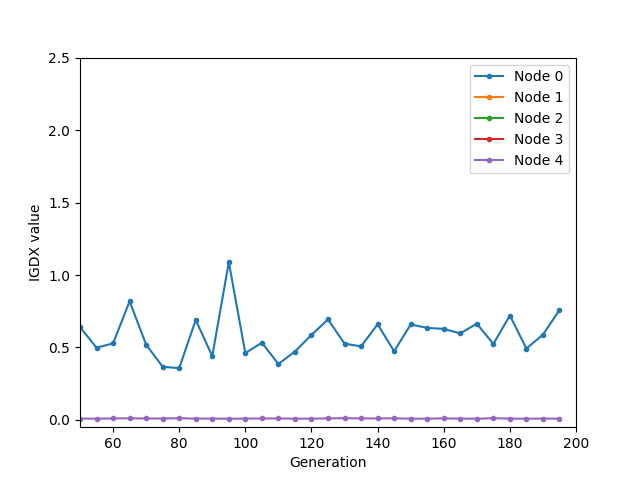}
        \caption{MOEA/D}
    \end{subfigure}
    \begin{subfigure}[b]{0.18\columnwidth}
        \centering
        \includegraphics[width = \columnwidth]{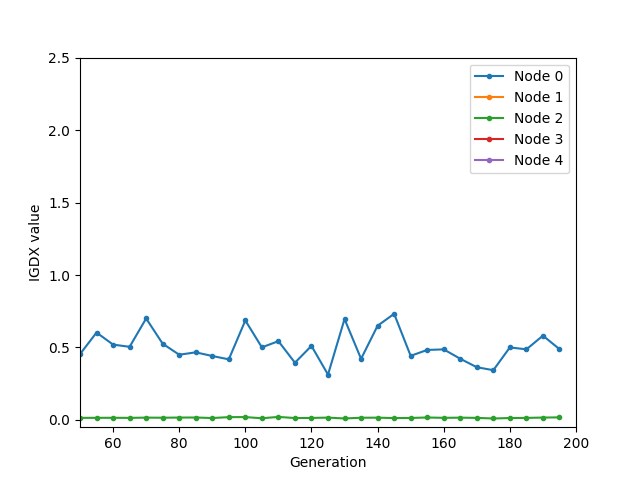}
        \caption{NSGA-II}
    \end{subfigure}
    \begin{subfigure}[b]{0.18\columnwidth}
        \centering
        \includegraphics[width = \columnwidth]{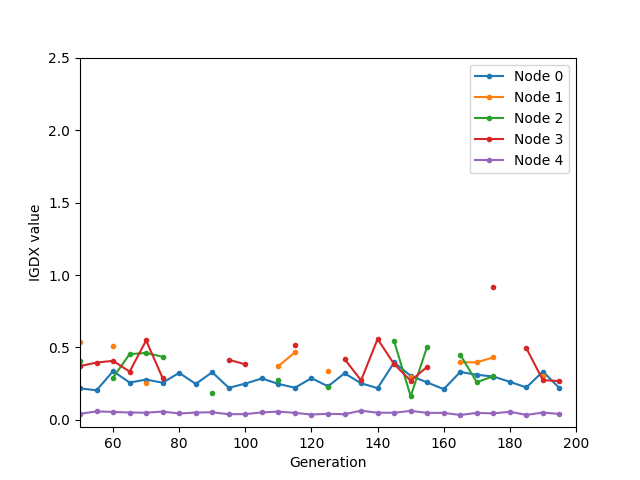}
        \caption{OMOPSO}
    \end{subfigure}
    \caption{Basin-wise IGDX for each basin of attraction (for a breadth-base basin graph)}
    \label{fig:basin-wiseIgdx for another example}
\end{figure}

For another example, we take a basin graph illustrated in \cref{fig: basin graph of breadth 1}, which we call a breadth-base basin graph. 
This 3BC has five nodes $\emptyset$, $1^+$, $1^-$, $2^+$, and $2^-$. 
As with the depth-base 3BC, label these nodes in order from 0 to 4. 
For this 3BC, the line graphs of the basin-wise IGDX for each node are shown in \cref{fig:basin-wiseIgdx for another example}. 

Comparing the graphs for the depth-base and breadth-base 3BCs in \cref{fig:basin-wiseIgdx} and \cref{fig:basin-wiseIgdx for another example}, we can find the following differences in performance of each solver for depth-base and breadth-base 3BCs, which implies 3BCs for different basin graphs have different properties as multimodal multi-objective benchmark problems.

\begin{enumerate}
    \item 
    GDE can only find node 1 for the depth-base 3BC, but can find the node 4 (global optima) for the breadth-base 3BC in around 130th generation. 
    \item 
    IBEA cannot find any node other than node 0 for the depth-base 3BC, but can find the global optima for the breadth-base 3BC. 
    \item 
    MOEA/D can find global optima and also find other local optima for the depth-base 3BC, and almost all generation points are located in the basin of the node 4 for the depth-based 3BC. 
    \item 
    NSGA-II cannot find (local) optima by the 200th generation for the depth-base 3BC, but is trapped in the basin for node 2 for the breadth-base 3BC. 
    \item 
    OMOPSO cannot find any node other than node 0 for the depth-base 3BC, but can find the global optima for the breadth-base 3BC. 
\end{enumerate}

\section{Conclusion}
We have proposed a novel benchmark problem suite called the 3BC, which has an explicitly definable multimodal landscape with an arbitrary dimensional design space. Our work employs the basin graph for generating the problem. By utilizing our benchmark problem, one can now define the positions as well as the value of the optima for benchmarking a multi-objective optimization solver. Additionally, our benchmark problem can be set to be non-separable.

\bigskip

\appendix

\section{On separability and shape of the Pareto front}
\label{sec: appendix dealing with separability}

\paragraph{Separability}

Both separable and non-separable problems can be obtained by our method.
A solver can detect a Pareto optimum in a separable problem just by optimizing each parameter individually \cite{huband2005scalable}. 
Hence, non-separability is one of the properties required for complex benchmark problems. 
Our benchmark is non-separable if $m_\emptyset$ is not the minimum among $m_s$'s proved in the following.

The definition of separability \cite{huband2006review} is the following. 
Let $O$ be a single objective optimization problem with design valuables $x_i$'s. 
For each index $i$ we can consider a single objective and single variable problem $P_{O, \textbf{x}, i}$ from $O$ by fixing a parameter $\textbf{x}$ of the design space and varying only $x_i$. 
Then $x_i$ is said to be \textit{separable} on $O$ if the set of the global optima $P^*_{O, \textbf{x}, i}$ of $P_{O, \textbf{x}, i}$ is the same for any parameter $x$. 
The optimization problem $O$ is said to be \textit{separable} if every $x_i$ is separable on $O$. 
Morerover, a multi-objective optimization problem is said to be \textit{separable} if every objective of the problem is separable.

One necessary condition for our benchmark problem being separable is that $m_\emptyset$ is the minimum along all $m_s$'s. 
$\tilde f$ is defined by \eqref{eq: def of rotated objective function}. 
\begin{equation*}
    \tilde f(t, x) \coloneqq 
    \frac{1}{\sqrt 2}
    \begin{pmatrix}
    1 & 1 \\ -1 & 1
    \end{pmatrix} 
    \begin{pmatrix}
    t \\ f(t, x)
    \end{pmatrix}.
    \tag{\ref{eq: def of rotated objective function}} 
\end{equation*}
Since $\tilde f$ is separable, the first component $\frac{1}{\sqrt{2}}(t + f(t,x))$ is separable. That each parameter $x_i$ is separable on this component means that $x_i$ is separable on $f$. 
For each $i$, $P_{f, (1, x), i} = |x_i| + \text{constant}$, where $x = (1/4, 0, \dots, 0)$ if $i\neq 1$ and $x = (0, 1/4, 0, \dots, 0)$ if $i = 1$. 
Hence $P^*_{f, (1, x), i} = \{0\}$. 
The separability of $f$ implies that 
\begin{equation*}\label{eq: separability of f if m_emptyset is the minimum}
    P^*_{f, (t, x), i} = \{0\}
\end{equation*}
for any $(t, x)$ and $i$. 
Now we prove that $m_\emptyset$ is the minimum.
Take $t_0> 0$ large enough for $f(t_0, x_s) = m_s$ for every node $s$. 
Let $x_s = (y_1, \dots, y_n)\in \mathbb R^n$ be the coordinates of $x_s$. 
By using \eqref{eq: separability of f if m_emptyset is the minimum} $n$ times, we have 
\begin{align*}
    m_s 
    &= f(t_0, y_1,\dots, y_n) \\
    &\ge f(t_0, y_1, \dots, y_{n-1}, 0) \\
    &\ge f(t_0, y_1, \dots, y_{n-2}, 0, 0) \\ 
    &\dots \\ 
    &\ge f(t_0, 0, \dots, 0) \\ 
    &= m_\emptyset. 
\end{align*}
Hence $m_\emptyset$ is the minimum. 




\paragraph{Concave-ity of the Pareto optimal shape}

The shape of a Pareto front often affects the performance of solvers. 
Each local Pareto front of our benchmark is convex (illustrated in \cref{fig: condition of Pareto front}). 
The Pareto front, which is the combination of local Pareto fronts, is generally a mix of linear, convex, and concave fronts. 

\bibliographystyle{amsplain}
\bibliography{main-en.bib}

\end{document}